\documentclass[10pt,twoside]{amsart}

\usepackage{amssymb}

\usepackage[english,francais]{babel}
\usepackage{amsfonts}
\usepackage{amssymb}
\usepackage{amsmath}
\usepackage{latexsym}
\usepackage{amsbsy}
\usepackage[latin1]{inputenc}
\usepackage[T1]{fontenc}
\newtheorem{e-proposition}[theorem]{Proposition}

\newtheorem{e-definition}[theorem]{Definition\rm}

\newtheorem{theoreme}{Th\'eor\`eme}[section]

\newtheorem{proposition}[theoreme]{Proposition}
\newtheorem{corollaire}[theoreme]{Corollaire}
\newtheorem{definition}[theoreme]{D\'efinition\rm}

\newcommand{\eps}{\varepsilon}

\newcommand{\Om}{\Omega}
\newcommand{\Vol}{\mathrm{Vol}}

\newcommand{\N}{\mathbb N}

\newcommand{\R}{\mathbb R}
\newcommand{\C}{\mathbb C}

\def\og{\leavevmode\raise.3ex\hbox{$\scriptscriptstyle\langle\!\langle$~}}
\def\fg{\leavevmode\raise.3ex\hbox{~$\!\scriptscriptstyle\,\rangle\!\rangle$}}

\begin{document}

%\selectlanguage{francais}

\title[]{Entre analyse complexe et superanalyse.}

% utiliser les étiquettes pour indiquer l'adresse de chaque auteur,
%     s'il y a plusieurs adresses

% \author[label1,label2]{}
% \address[label1]{}
% \address[label2]{}

\maketitle

\begin{center}

\author{Pierre BONNEAU, Anne  CUMENGE}

\end{center}

\vspace{0,5cm}

\selectlanguage{francais}

\begin{abstract}
\selectlanguage{francais} Dans le cadre de la superanalyse, nous
obtenons une théorie des
 fonctions voisine de l'analyse complexe dès lors que nous imposons
aux superalgèbres réelles considérées une condition (A) qui
généralise en dimension supérieure à deux la relation $1+i^{2}=0$
de $\C$. la condition (A) permet l'obtention d'une formule de
représentation intégrale pour les fonctions superdifférentiables.
Nous donnons entre autres un théorème de super-différentiabilité
séparée et un théorème de prolongement de type Hartogs-Bochner
pour les fonctions superdifférentiables.

\end{abstract}

\vskip 0.5\baselineskip

\selectlanguage{english}
\begin{abstract}

%\noindent{\bf Abstract}
\vskip
0.5\baselineskip
{\bf Midway between complex analysis and
superanalysis. } In the framework of superanalysis we get a
functions theory close to complex analysis,
 under a suitable condition (A) on the real superalgebras in consideration
 (this condition  is a generalization of the classical relation $1+i^{2}=0$ in
 $\C$). Under the condition (A), we get an integral
 representation formula for the superdifferentiable functions.
 We give a result of Hartogs type of separated superdifferentiabity and
and a continuation theorem  of Hartogs-Bochner type for the superdifferentiable functions.
\end{abstract}

\bigbreak\noindent {Rubriques : Analyse complexe/Analyse sur les
superespaces et les espaces gradués.}

\vspace{1cm}

%%%%%%%%%%%%%%%%%%%%%%%%%%%%%%%%%%%%%%%%%%%%%%%%%%%%%%%%%%%%%%
\selectlanguage{francais}
% texte principale
%\section{Introduction}
%\label{}
\section{Introduction et notations}
Conformément à \cite{K}, nous noterons $\Lambda$ une superalgèbre
 commutative (en abrégé CSA), c'est-à-dire un $\mathbb{R}$-espace
vectoriel $\Lambda=\Lambda_{0}\bigoplus\Lambda_{1},$  somme
directe de deux sous-espaces vectoriels $\Lambda_{0}$ et
$\Lambda_{1}$,  et muni d'une structure d'algèbre associative
unitaire telle que $\Lambda_{0}$ est une sous-algèbre de
$\Lambda$, le produit de deux éléments de $\Lambda_{1}$ est dans
$\Lambda_{0}$, le produit d'un élément de $\Lambda_{0}$ par un
élément de $\Lambda_{1}$ est dans $\Lambda_{1}$, un élément de
$\Lambda_{0}$ commute avec tout élément de $\Lambda$ alors que
deux éléments de $\Lambda_{1}$ anti-commutent. Nous définissons le
$\Lambda_{1}$-annihilateur comme étant
$^{\bot}\Lambda_{1}=\{\lambda\in\Lambda :\lambda\Lambda_{1}=0\}.$
Toutes les CSA considérées ici seront de dimensions finies. Nous
noterons $(e_{0}, e_{1}, ...,e_{p})$ une base de $\Lambda_{0}$
($e_{0}$ est l'élément unité de $\Lambda$), $(\eps_{1},\eps_{2},
...,\eps_{q})=(e_{p+1},e_{p+2}, ...,e_{p+q})$ une base de
$\Lambda_{1}$, et définissons les coefficients de structure
(d'algèbre) $\Gamma$ de $\Lambda$ par
$e_{i}e_{j}=\sum_{k=0}^{p+q}\,\Gamma_{i,j}^{k}e_{k}$ pour
$i,j=0,1, ...,p+q$.

\smallbreak On appelle superespace sur la CSA $\Lambda$ tout
$\mathbb{R}$-espace vectoriel $\mathbb{R}_{\Lambda}^{n,m}\,=
\Lambda_{0}^{n}\times\Lambda_{1}^{m}$
 où $n,m\in\mathbb{N}$. \\
\begin{definition}: Soit $U$ un ouvert de
$\mathbb{R}_{\Lambda}^{n,m}$ ; une application  $F$   de $U$ dans
$\Lambda$ {\it est superdifférentiable à droite} ( {\it
S-différentiable} en abrégé) en $x\in U$ s'il existe des éléments
$\frac{\partial F}{\partial x_{j}}(x)$ de $\Lambda$, $j=1,...,n+m$
tels que, pour tout $h\in\mathbb{R}_{\Lambda}^{n,m}$ tel que
$x+h\in U$, on ait $F(x+h)=F(x)\,+\,\sum_{j=1}^{n+m}\frac{\partial
F}{\partial x_{j}}(x)h_{j}\,+\,o(h)$ avec $lim_{\|h\|\mapsto
0}\frac{\|o(h)\|}{\|h\|}=0$ où $\|\,.\,\|$ est une norme sur
l'espace  $\mathbb{R}_{\Lambda}^{n,m}.$
\end{definition}

\smallbreak\noindent Nous soulignons que $\frac{\partial
F}{\partial x_{1}}(x),...,\frac{\partial F}{\partial x_{n}}(x)$
sont définis de façon unique  tandis que $\frac{\partial
F}{\partial x_{n+1}}(x),...,\frac{\partial F}{\partial
x_{n+m}}(x)$ le sont seulement modulo $^{\bot}\Lambda_{1}.$
\\
La condition de S-différentiabilité de $F$ traduit que $F$ est
Fréchet différentiable avec une différentielle définie par des
opérateurs de multiplication par des éléments de $\Lambda$.

\medbreak Les deux hypothèses suivantes $\bf{(A_{0})}$ et
$\bf{(A_{1})}$ interviendront de manière naturelle dans la
recherche d'un noyau reproduisant pour les fonctions
S-différentiables.

\smallbreak\noindent $\bf{(A_{0})}$:  il existe une base
$(e_{0}=1,e_{1},...,e_{p})$ de $\Lambda_{0}$ vérifiant
$\sum_{k=0}^{p}\,e_{k}^{\;2}=0$

 \smallbreak\noindent $\bf{(A_{1})}$:  il
existe une base $(\eps_{1},...,\eps_{q})$ de $\Lambda_{1}$ et une
suite finie $s_{1}=1 <s_{2}<...<s_{r}<s_{r+1}=q+1$ telles que,
pour tout $j=1,...,q$, il existe $a_{j}\in\Lambda_{0}$ vérifiant
$\eps_{j}=a_{j}\eps_{s_{k}}\;si\;s_{k}\leq j<s_{k+1}$, avec
$a_{s_{1}}=\dots =a_{s_{r}}=e_{0}$ et
$\sum_{j=s_{k}}^{s_{k+1}-1}\,a_{j}^{\;2}=0\;\mbox{ pour tout }
k=1,2,...,r.$

\medbreak\noindent On écrit  un élément $x$ de
$\Lambda_{0}^{n}\times\Lambda_{1}^{m}$ :
  $x=(y,\theta)=(y_{1},...,y_{n},\theta_{1},...,\theta_{m}),
\mbox{ où } y_{i}=\displaystyle{\sum_{j=0}^{p}\,y_{i}^{j}e_{j}
\mbox{ et } \theta_{k}=\sum_{l=1}^{q}}\,\theta_{k}^{l}\eps_{l}\,.$

\begin{definition}:  L'opérateur de Cauchy-Riemann $d''$  sur $\Lambda_0^{\, n}\times \Lambda_1^{\, m}$ est
défini par :
$$d''=\sum_{i=1}^{n}\sum_{j=1}^{p}\,dy_{i}^{j}
\left(\frac{\partial }{\partial y_{i}^{j}}-e_{j}\frac{\partial
}{\partial y_{i}^{0}}\right)\,
+\,\sum_{l=1}^{m}\sum_{k=0}^{r}\sum^{s_{k+1}-1}_{i=s_{k}+1}\,d\theta_{l}^{i}
\left(\frac{\partial
}{\partial\theta_{l}^{i}}-a_{i}\frac{\partial}{\partial
\theta_{l}^{s_{k}}}\right)$$
\end{definition}

Dans le cas le plus simple, lorsque $m=0$ et $n=1$, nous avons
cherché à définir sur $\Lambda_0$ un opérateur $d''$ de type
Cauchy-Riemann dont le noyau coïncide, dans sa partie constituée
de 0-formes, avec les fonctions S-différentiables. Contrairement à
ce qui se passe en analyse complexe, nous n'avons pas, faute de
conjugaison, de moyen canonique de choisir $d''$ ; la définition
que nous avons retenue ne correspond pas au choix effectué dans
$\mathbb{C}$ pour $\bar\partial$, mais nous avons bien $\ker d'' =
\ker \bar\partial$ dans ce cas.

\smallbreak \noindent {\sc Remarque}: si  $F$ est S-différentiable
sur  $U$ ouvert alors $d''F=0$ sur $U$ avec équivalence si
$U\subset \Lambda_0^{\, n}$.

\smallbreak
\begin{definition}: Soit $U$  un ouvert de
$\mathbb{R}_{\Lambda}^{n,m}$ ; une application $F$  de $U$ dans
$\Lambda$ est quasiment S-différentiable (ou {\it
qS-différentiable}) sur $U$ si elle est Fréchet-différentiable sur
$U$ et vérifie $d''F=0$.
\end{definition}

\smallbreak \noindent {\sc Exemples }: {(1) } Si $m=0$ et
$\Lambda_0=\C$, la superanalyse est l'analyse complexe ; $f$ est
S-différentiable si et seulement  $f$ est $\C$-analytique, donc
très
régulière.\\
{(2)} Si $m=0$ et $n=1$ avec $\Lambda_0 =Vect (e_0 ,e_1 )$ où $e_0
=1$ et $e_1^{\,2}=e_0$, la superanalyse est l'analyse hyperbolique
; dans ce cas, il se peut qu'une
 fonction S-différentiable ne soit pas plus régulière que
  $C^1$ au sens de Fréchet.\\
{(3)} Soit $A$  une algèbre de Clifford réelle de dimension $2^k$
dont les  générateurs   $e_1 ,..., e_k $  vérifient $e_i e_j + e_j
e_i = -2 \delta_{ij}$ ; $A$ peut s'écrire comme une superalgèbre
$A=\Lambda_0 \oplus \Lambda_1$ avec
 $\Lambda_0 =Vect
\{e_I \; ;\, |I| \mbox{ pair } \}$, $\Lambda_1 = Vect \{e_I\; ;\,
|I| \mbox{ impair }\}$, où $e_\emptyset =1$,et si $I=(i_1 , \dots
, i_\nu )\;,\; e_I =e_{i_1 } \dots e_{i_\nu} $.

\smallbreak\noindent  La condition ($A_0$) est vérifiée dans
l'exemple 1 et, lorsque $k \equiv 2$ modulo 4, dans l'exemple 3.
Elle ne l'est pas dans le second (ce serait d'ailleurs en
contradiction avec la proposition \ref{harmonicité} du paragraphe
suivant).

\medbreak {\it Dans toute la suite, nous supposerons  satisfaites
les deux hypothèses $\bf{(A_{0})}$ et $\bf{(A_{1})}$.}

\section {Représentation intégrale de formes et fonctions.}
 \begin{theoreme}
L'opérateur de Cauchy-Riermann $d"$ dans
$\mathbb{R}_{\Lambda}^{n,m}$ admet une solution fondamentale
$\Omega$  donnée par:
\begin{eqnarray}
\Omega=\frac{c(n,m,p,q)}{\|x\|^{n(p+1)+qm} }  \Big[\sum_{i=1}^{n}\sum_{j=1}^{p}\,(-1)^{(p+1)(i-1)+j}(y_{i}^{j}e_{0}
+y_{i}^{0} e_{j})dy_{1}^{0}...\widehat{dy_{i}^{j} }...dy_{n}^{p}d\theta^{1}_{1}...d\theta^{q}_{m}\, \nonumber \\
+\,\sum_{l=1}^{m}\sum_{k=1}^{r}\sum^{s_{k+1}-1}_{i=s_{k}+1}\,(-1)^{n(p+1)+i-1}
(\theta_{l}^{i} e_{0}+\theta_{l}^{s_{k}} a_{i})dy^{0}_{1}.
..dy^{p}_{n}d\theta_{1}^{1}...
\widehat{d\theta_{l}^{i}}...d\theta_{m}^{q} \Big]
\end{eqnarray}
\noindent où $c(n,m,p,q)= -\big(( n(p+1)+qm ) \Vol
(B(0,1))\big)^{-1}$
\end{theoreme}

\noindent {\sc Schéma de preuve }: dans le cas où
$\mathbb{R}_{\Lambda}^{n,m}=\Lambda_{0}$,  nous n'avons pas de
conjugaison comme dans $\mathbb{C}$, et pas non plus l'intégrité
de $\Lambda_{0}\neq \C$ ; nous cherchons alors une solution
fondamentale de l'opérateur $d''$ de la forme
$\Omega(x)=\frac{A(x)}{\|x\|^{p+1}}$ où $A(x)$ est une p-forme à
coefficients polynômes homogènes de degré 1. Nous sommes alors
conduits, de façon naturelle, à imposer la condition ($A_{0}$) à
$\Lambda_{0}$ afin d'obtenir une solution fondamentale explicite
de l'opérateur $d''$. Le cas $\Lambda_0^{\, m}$
est seulement plus technique.\\
Si $\mathbb{R}_{\Lambda}^{n,m}=\Lambda_{1}$, n'ayant pas, dans
$\Lambda_{1}$, d'élément unité (ni de conjugaison comme dans le
cas des quaternions) nous sommes amenés, pour définir un opérateur
de Cauchy-Riemann contenant, dans son noyau, les fonctions
S-différentiables, à particulariser certains éléments et sommes
conduits à la condition ($A_{1}$). Un calcul direct fournit alors
une solution fondamentale dans $\Lambda_{1}$, puis dans
$\Lambda^{m}_{1}$, et enfin dans $\mathbb{R}_{\Lambda}^{n,m}$.

\medbreak Soit D un ouvert de $\mathbb{R}_{\Lambda}^{n,m}$ borné à
frontière lisse et
$\Psi:\,\overline{D}\times\overline{D}\rightarrow\mathbb{R}_{\Lambda}^{n,m}$
définie par $\Psi(x',x)=x'-x$. Définissons $K(x',x)={\Psi^{\ast}}
\Omega$.
 Le noyau $K$ vérifie, si
l'on pose $d''=d_{x'}''+d_{x}''$:\\
$d''K(x',x)=d''\Psi^{\ast}\Omega=\Psi^{\ast}d''\Omega={\Psi^{\ast}}[0]=[\Delta]$
où $[\Delta]$ désigne le courant d'intégration sur la diagonale.

\smallbreak Nous obtenons alors classiquement (cf. \cite{HP} par
exemple ) :

\begin{corollaire}\label{corol-repres}
   Si f est une forme (ou une
fonction) continue sur $\overline{D}$ et de classe $C^{1}$ dans D,
alors pour tout $x\in D$, nous avons :
\begin{eqnarray}\label{repres-int}
f(x)=\int_{\partial D}\,f(x')K(x',x)\,-\,\int_{D}\,d''f(x')K(x',x)\,+\,d_{x}''\int_{D}\,f(x')K(x',x).
\end{eqnarray}
\end{corollaire}

Dans la formule (\ref{repres-int}) le dernier terme est nul si f
est une fonction, les deux derniers le sont si f est une fonction
qS-différentiable.

\medbreak Il  découle de la représentation intégrale obtenue pour
les fonctions S-différentiables des propriétés voisines de celles
des fonctions holomorphes, propriétés que nous précisons dans les
deux derniers paragraphes.

\section{ Propriétés d'analyticité des fonctions S-différentiables}

 \begin{proposition}\label{harmonicité}
 Si $f:\,D\subset \Lambda_{0}^{n}\times\Lambda_{1}^{m}\rightarrow\Lambda$ est
 qS-différentiable sur $D$, alors elle est harmonique sur $D$, donc $\mathbb{R}$-analytique sur $D$.
\end{proposition}

\noindent{\sc notations}. Si $1\leq j\leq m$ et $\theta_j =
\displaystyle{\sum_1^{q} \theta_j^i \eps_i}$ :\;\\
$Z(\theta_{j}):=\displaystyle{\sum_{i=1}^{q}\,a_{i}\theta_{j}^{i}}$
\; et $\pi_k (\theta_j ):=\displaystyle{\sum_{i=s_k}^{s_{k+1}
-1}}\theta_j^i\, \eps_i
 \,,\,\forall k=1,\dots ,r$

\smallbreak  Une fonction qS-différentiable $f$ est analytique
réelle en les variables $(y^0_{1},..., y^p_{1},\dots,\\
 y^0_{n},\dots ,  y^p_{n} ,\theta^1_{1},\dots,$ $\theta^q_1 ,\dots
, \theta^1_{m} ,\dots , \theta^q_m)$, mais possède en fait une
propriété d'analyticité bien plus forte :

\smallbreak
\begin{proposition}\label{S-analyticite}
Soit
$f:\,D\subset\Lambda_{0}^{n}\times\Lambda_{1}^{m}\rightarrow\Lambda$
une fonction qS-différentiable sur $D$. Pour tout
$(b,\beta)=\\
(b_{1},...,b_{n},\beta_{1},...,\beta_{m})\in D$, il
existe $R>0$ et des scalaires $A_{I,J}\in\Lambda$ où I et J sont
des multi-indices dans $\mathbb{N}^{n}$ et $\mathbb{N}^{m}$
respectivement, tels que
pour $\|y_{j}-b_{j}\|<R,\;j=1,...,n;\;\;\|\theta_{j}-\beta_{j}\|<R,\;j=1,...,m$\\
$$f(y,\theta)\,=\,\sum_{I,J_{1} ,\dots , J_{r}}\,A_{I,J}(y-b)^{I}    (Z_1  (\theta -\beta ))^{J_{1}}...(Z_r (\theta -\beta ))^{J_{r}}   $$\\
avec absolue convergence de la série.
\end{proposition}

\smallbreak\noindent Nous avons noté
$(y-b)^{I}=(y_{1}-b_{1})^{i_{1}}...(y_{n}-b_{n})^{i_{n}}\mbox{ si }\,I=(i_{1},...,i_{n})\in\mathbb{N}^{n}$\\
$Z_k  (\theta -\beta ) = \big( Z(\pi_k  ( \theta_{1}-\beta_{1})), \dots , Z(\pi_k ( \theta_{m}-\beta_{m})) \big) $ \\
et $Z_k  (\theta -\beta ) ^{J_{k}} = Z(\pi_k
(\theta_{1}-\beta_{1}))^{\, j^k_1} \dots Z(\pi_k (
\theta_{m}-\beta_{m}))^{\, j^k _{m}}$ si $J_k =(j^{k} _1, \dots ,
j^{k}_{m} )\in \N^m$.

\section {Deux théorèmes de type  Hartogs. }

%Un théorème de S-différentiabilité séparée.

\smallbreak Une fonction $f$ séparément holomorphe par rapport à
chaque variable $z_j \,$ de $\C^n$ ,\,$ j=1,\dots , n$,  sur un
domaine $\Om$ de $\C^n$ est holomorphe sur $\Om$ sans autre
hypothèse de régularité globale sur $f$ mais ce théorème de
Hartogs n'est plus
valable pour les fonctions  $\R$-analytiques. \\
%Pour les fonctions S-différentiables,
Lorsque les coefficients $\Gamma_{ij}^k$ vérifient une condition
($\mathcal{P}$) de positivité, nous obtenons un résultat de type
Hartogs de S-différentiabilité séparée.

\smallbreak
\begin{definition} une CSA $\Lambda_0$ vérifie la condition de
positivité ($\mathcal{P}$) si, pour tous $X,Y\in \R^{p+q+1}$, on a
\\
 $ \big(   \sum_0^{p+q} (X_j )^2 \big)
\sum_{\underset{j=0,..,p+q}{ k=0,...,p}} \big(\sum_{m=0}^{p+q} Y_m
\Gamma_{mk}^j\big)^2
 -2 \sum_{k=0}^p\, \big( \sum_j X_j \sum_m Y_m \Gamma_{mk}^j \big)^2 \,
\geq 0 \,.$
\end{definition}

\smallbreak
\begin{theoreme}\label{hartogs-separe}
On suppose que la CSA $\Lambda$ vérifie la propriété
($\mathcal{P}$) et toujours la condition $(A_0 )$. Alors une
fonction $f$ définie sur un domaine $D$ de $\Lambda_0^{\,n}$ à
valeurs dans $\Lambda$ séparément S-différentiable en les
variables $y_1 , \dots , y_n$ est de classe $C^1$ au sens de
Fréchet sur $D$ et par suite est S-différentiable en tout point
$y=(y_1 ,\dots ,y_n )$ de $D$.
\end{theoreme}

 \noindent {\sc remarque} :
  nous utilisons la sous-harmonicité de $\log ||f||$, qui
 découle de la condition ($\mathcal{P}$),
 pour en déduire la continuité de
 la fonction $f$. Or cette propriété de continuité ne dépend
 pas de la norme choisie, ni de la métrique retenue
 pour définir le laplacien. Nous pouvons ainsi remplacer la condition ($\mathcal{P}$) par la condition
 suffisante ($\mathbb{P}$) suivante :
il existe une base $(e'_{k})$ de $\Lambda_{0}$ vérifiant
$\sum^{p}_{k=0}\,(e'_{k})^{2}=0$, et  une base $(e''_{\alpha})$ de $\Lambda$
telles que, pour tout
$X=\sum^{p+q}_{\alpha=0}\,X_{\alpha}e_{\alpha}$
et Y de $\Lambda$,\\
$\|X\|^{2}\sum_{k,\alpha}\big(\left\langle e'_{k}Y,
e''_{\alpha})\right\rangle \big)^{2}
\,-\,2\sum_{k}\big(\sum_{\alpha}\left\langle e'_{k}Y,
e''_{\alpha}\right\rangle X_{\alpha}\big)^{2}\;\geq\;0.$
\\

\begin{theoreme}: Sous les conditions ($A_0$) et ($A_1$), si $\partial \Om$ est le bord connexe d'un domaine
$\Om$ borné de $\Lambda_0 ^{\,n} \times \Lambda_1 ^{\, m}$ et $f$une fonction
qS-différentiable au voisinage de $\partial \Om$, alors $f$ se
prolonge en une fonction qS-différentiable sur $\Om$.
\end{theoreme}

\vspace{2cm} \bigbreak\noindent
%\author{}\\
\address{Equipe Emile Picard, Institut de Mathématiques, {\sc umr} 5219\\
 Universit\'{e} Paul Sabatier - 31062 Toulouse Cedex 9}\\
\email{{\it pierre.bonneau@math.univ-toulouse.fr}, {\it
anne.cumenge@math.univ-toulouse.fr}}

\end{document}